\documentclass[12pt,twoside]{article}
\usepackage{amsmath}
\usepackage{amssymb}
\usepackage{enumerate,epsfig}
\usepackage{epstopdf}
\usepackage{graphicx}
\newtheorem{lemma}{Lemma}[section]
\newtheorem{theorem}{Theorem}[section]
\newtheorem{corollary}[lemma]{Corollary}
\newtheorem{remark}{Remark}[section]
\newtheorem{prop}{Proposition}[section]
\newtheorem{defi}{Definition}[section]

\newcommand{\R}{ {\mathbb R} }
\newcommand{\N}{{\mathbb N}}
\newcommand{\eps}{{\varepsilon}}
\newcommand{\E}{{\mathbb{E}}}
\renewcommand{\div}{\mbox{div}\,}
\makeatletter
\@addtoreset{equation}{section}
\makeatother

\newcommand{\cqfd}{{\unskip\kern 6pt\penalty 500
\raise -2pt\hbox{\vrule\vbox to 6pt{\hrule width 6pt
\vfill\hrule}\vrule}\par}}

\newcommand{\ind}{{\mathbb I}}
\title{Critical non Sobolev regularity for continuity equations with rough force fields}
\author{
Pierre--Emmanuel Jabin
\footnote{CSCAMM and Dept. of Mathematics, University of Maryland,
College Park, MD 20742, USA. P.--E.~{\sc Jabin} is partially supported by NSF Grant 1312142 and by NSF Grant RNMS (Ki-Net) 1107444.}
}
\date{}
\begin{document}
\maketitle

\begin{abstract}
We present a divergence free vector field in the Sobolev space $H^1$ such that the flow associated to the field does not belong to any Sobolev space. 
The vector field is deterministic but constructed as the realization of a random field combining simple elements. This construction illustrates the optimality of recent quantitative regularity estimates as it gives a straightforward example of a well-posed flow which  has nevertheless only very weak regularity.
\end{abstract}

\noindent{\bf MSC subject classifications:} 34C11, 35L45, 37C10, 34A12, 34A36.
\section{Introduction}
Given a vector field $u(t,x):\;\R_+\times\R^d$, we consider the flow $X(t,s,x)$ associated to it
\begin{equation}
\partial_t X(t,s,x)=u(t,X(t,s,x)),\quad X(t=s,s,x)=x.\label{theflow}
\end{equation}
The well-known Cauchy-Lipschitz theorem implies that if $u$ is locally Lipschitz, then Eq. \eqref{theflow} is locally well posed for every choice of $x$.  The flow $X$ is in that case also locally Lipschitz in $x$ and this regularity is intimately connected with the existence and uniqueness result.

Well posedness for {\em non Lipschitz} velocity fields $u$  was first obtained by {R.J. Di Perna} and P.--L. { Lions} in \cite{DL} with the introduction of renormalized solutions for the associated continuity equation
\begin{equation}
\partial_t \rho+\div(\rho\,u)=0,\quad \rho(t=0,x)=\rho^0(x).\label{continuity}
\end{equation}
The result \cite{DL} required $u\in L^1_t W^{1,1}_x$ with appropriate bounds on $\div u$. 

In the general case of \eqref{theflow}, well posedness was proved with only $u\in BV$ and $\div u\in L^p$ in the seminal article \cite{Am} (see \cite{Bo} for a first result in the kinetic context). The assumption on $\div u$ can be replaced by lower and upper bound on a solution to \eqref{continuity} as shown in \cite{ADM}, though one then typically requires $u\in SBV$ (see also \cite{Ja}).

The optimality of $u\in BV$ for such a general setting was demonstrated in \cite{DeP} with an example of almost $BV$ vector field leading to non unique flows. We refer the readers to the two excellent summaries in \cite{DeL} and more recently \cite{AmCr}.

Contrary to the classical Cauchy-Lipschitz theory, well posedness results that rely on the theory of renormalized solutions do not in general provide  explicit regularity estimates on the flow $X$. In contrast a purely quantitative approach was introduced in \cite{CD}, focusing on estimates like
\begin{equation}
\int_{\R^d} \oint_{B(x,h)}\log\left(1+\frac{|X(1,0,x)-X(1,0,y)|}{h}\right)\,dy\,dx
\leq C_d\,\|u\|_{L^1_t W^{1,p}_x}.\label{optreg}
\end{equation}
Such direct estimates on the flow give both quantitative regularity results and well posedness (see also \cite{HaLL,LL,LL2} for another approach of well posedness working directly with the flow). 

A regularity result like \eqref{optreg} also implies some regularity on the solution to the continuity equation \eqref{continuity}. Regularity estimates have been obtained directly from Eq. \eqref{continuity} in \cite{BeJa} and in the more general context of non linear continuity equations but they require more complicated tools.

The regularity given by \eqref{optreg} is quite weak, a sort of $\log$ of a derivative and very far from the Lipschitz regularity obtained if $u\in W^{1,\infty}_{loc}$. It is therefore a very natural question whether this regularity is optimal or if something better should be expected. 

It is relatively easy to construct examples of flows with weak regularity with just the requirement that $u\in L^2_t H^1_x$ but no assumption on $\div u$. 
But there are very few examples if one imposes that $\div u\in L^\infty$ or even that {\em $u$ be divergence free}. To the author's knowledge, the only such other example was obtained in \cite{ACM}; other interesting examples of velocity fields occur in the slightly different but related context of mixing (see below).

Note that the assumption $\div u=0$ is an especially strong constraint in low dimensions: Obviously in dimension $1$ it implies that $u$ is constant and then $\rho(t,.)$ is as smooth as $\rho^0$ so that \eqref{optimalreg} cannot stand. 

This article presents such an example of weak regularity in the minimum dimension where it is possible, namely $d=2$.
\begin{theorem}
For any $\gamma>0$, there exists a divergence free, compactly supported velocity field $u\in L^2(\R_+,\ H^1(\R^2))\cap L^\infty$ s.t. the flow solving
\[
\partial_t X(t,0,x)=u(t,X(t,0,x)),\quad X(t=0,0,x)=x,
\]
does not belong to any Sobolev space $W^{1,p}$, $p\geq \gamma$, and more precisely
\begin{equation}
\limsup_{h\rightarrow 0}\int_{\R^2} \oint_{B(x,h)}\frac{|X(1,0,x)-X(1,0,y)|^\gamma}{h^\gamma}\,dy\,dx=+\infty.
\label{optimalreg}\end{equation}
\label{counterexample}\end{theorem}
\begin{remark}
The choice of $u\in L^2_t H^1_x$ was just a matter of convenience and simplicity. The example could easily be adapted to $u\in W^{1,p}_x$ for any $p<\infty$. 
\end{remark}

\begin{remark}
Note that, taking $\gamma<1$, \eqref{optimalreg} shows that the flow $X$ cannot be in $BV$ nor in any generalized Sobolev spaces $W^{1,p}$ with $\gamma<p<1$. In fact from the construction, it can be checked that $X$ is Lipschitz regular in a set $\Omega$ with constant $\exp (C\,\Omega^c)$. This is exactly the regularity predicted by \cite{CD}.
\end{remark}

\begin{remark}
The construction of $u$ relies on the choice of one $\gamma>0$, though obviously one particular field $u$ would also work for any $\gamma'\geq \gamma$. Extending the construction so that \eqref{optimalreg} holds with the same $u$ for any $\gamma>0$ seems non trivial. 
\end{remark}

Though in the general setting, $u\in BV_x$ or $u\in W^{1,p}$ is critical, there are several cases where well posedness has been obtained with less regularity. 

In dimension $2$ for time independent fields $u$, well posedness holds with essentially no regularity on $u$. This was first noticed in \cite{BoDe} in the kinetic context for a continuous $u$, extended to $u\in L^p$ in \cite{Ha1}. It was later extended to  $u$ with only bounded divergence in \cite{CR} and \cite{CCR} with the most general result available in \cite{ABC}.
 
Structural assumptions on $u$ can also compensate for lack of regularity. Well posedness was for instance recently proved in \cite{BoCr} for fields $u$ which result from the singular integral of a measure. 

The kinetic setting also provides important additional structure: It was proved in \cite{CJ} that it is enough for well posedness in any dimension with only $u\in H^{3/4}$. This was later developed in \cite{JaMa} for $u\in L^2_x$ provided that it solves a wave equation (such as Maxwell's system).

In the cases mentioned above where the critical scale is a different space, it is not obvious what the optimal regularity for the flow should be or if better regularity could be obtained above the critical space but below $W^{1,\infty}_{loc}$.

\bigskip

Note as well that, from Theorem \ref{counterexample}, it is straightforward to deduce on solutions to the continuity equation
\begin{corollary}
For any $\gamma>0$, there exists a divergence free, compactly supported velocity field $u\in L^2(\R_+,\ H^1(\R^2))\cap L^\infty$ and an initial data $\rho^0\in W^{1,\infty}(\R^2)\cap L^1(\R^2)$ s.t. the unique solution to
\[
\partial_t \rho+\div(\rho\,u)=0,\quad \rho(t=0,x)=\rho^0(x),
\]
does not belong to any Sobolev space $W^{1,p}$, $p\geq \gamma$, and more precisely
\begin{equation}
\limsup_{h\rightarrow 0}\int_{\R^2} \oint_{B(x,h)}\frac{|\rho(1,x)-\rho(1,y)|^\gamma}{h^\gamma}\,dy\,dx=+\infty.
\label{optimalreg2}\end{equation}
\label{counterexample2}
\end{corollary}
%
As mentioned above the main point of comparison for Theorem \ref{counterexample} is the example obtained in \cite{ACM}. The construction behind Theorem \ref{counterexample} is random at its base. Roughly speaking, one builds a random field $u^{rand}$ s.t. the associated flow satisfies
\[ 
\E \int_{\R^2} \oint_{B(x,h)}\frac{|X(1,0,x)-X(1,0,y)|^\gamma}{h^\gamma}\,dy\,dx=f(h),
\]
for some function $f(h)\rightarrow \infty$ as $h\rightarrow 0$. Then one only has to choose one realization of $u^{rand}$. 

The advantage to this approach is that it relatively simple to build, using only elementary fields and transforms (rotations in fact). This is because the stochasticity  takes care of most of the difficult choices. Combining the various random angles in the construction for instance would be very difficult to do explicitly but is straightforward here as it is enough to explain what happens on average. 

This type of approach may also prove useful in the future to prove some sort of generic regularity. 

Its main drawback is that it is not explicit. Much about the flow remains in the dark: Its choice itself depends on the exponent $\gamma$ that is used for example. In particular it is not possible to have the same level of information and precision as in the construction in \cite{ACM}.

As a consequence, and contrary to \cite{ACM}, the example given here does not say much about {\em mixing}. In several respect, mixing is the dual question of the one investigated here. Given a class ${\cal U}$ of vector fields and a bounded and smooth initial data $\rho^0$, denote
\[
{\cal S}_{u,\rho^0}=\{\rho(t,x)\in L^\infty\ \mbox{solution to \eqref{continuity} with initial data}\ \rho^0\}.
\]
Then one typically wishes to estimate 
\[
M(t,\rho^0)=\inf_{u\in {\cal U}} \inf_{\rho\in{\cal S}_{u,\rho^0}} \|\rho(t,.)\|_{\dot H^{-1}_x}.
\]
This is an important problem in itself with many applications in Fluid Mechanics in particular. 

If ${\cal U}$ only includes divergence free vector fields then \eqref{continuity} propagates any $L^p$ norm of $\rho$ and in particular $\|\rho(t,.)\|_{L^2}$ is constant in time. By interpolation, this implies that $M(t,\rho^0)$ can only be small if $\rho$ is not too smooth, for instance if $\|\rho(t,.)\|_{H^1_x}$ is large.

Therefore regularity results for \eqref{continuity} imply lower bounds on $M(t,\rho^0)$. For instance the equivalent of the $\log $ estimate \eqref{optreg} of \cite{CD} is the lower bound $\|\rho(t,.)\|_{\dot H^{-1}_x}\geq \frac{e^{-Ct}}{C}$ proved in \cite{IKX,Seis}. 

Similarly perfect mixing $\|\rho(t,.)\|_{\dot H^{-1}_x}=0$ can only occur in finite time if $\rho$ loses all regularity. It was observed in \cite{LTD} and proved in \cite{LLNMD} for 
\[
{\cal U}=\{u\in L^\infty_t\,L^2_x,\quad \div u=0\}, 
\]
which corresponds to ill posedness for the PDE \eqref{continuity} and the flow \eqref{theflow}.

This connection between mixing and critical regularity for \eqref{continuity} 
was powerfully made in \cite{Bress}, together with important (and still open) conjectures. Examples of mixing, such as the recent \cite{YaZl} giving the optimal rates for $u\in W^{1,p}$, could therefore turn out to be very useful as well for the question of critical regularity.

Nevertheless note that  a good mixing solution $\rho$ could still in principle belong to Sobolev spaces $W^{s,p}$ (even with $s$ large); mixing only implies that the Sobolev norm has to be large. Conversely examples such as given here do not necessarily imply mixing (but the example in \cite{ACM} does). 
\section{Sketch of the construction and proof of Theorem \ref{counterexample}}
\subsection{Main steps of the construction}
The basic idea of the construction is to make sure that the trajectories will go through a certain basic random structure many times, and each time the distance will be amplified by a given factor $c>1$.

The result of the construction is summarized by
\begin{prop}
For any $a\in\R^2$, any $\eps>0$, $\alpha<1$ and any $\alpha<\beta<1$, there exists a divergence free vector field $u_{a,\eps}$ with support in the ball $B(a,\;2\,\eps)$ and the bound
\begin{equation}
\|u_{a,\eps}\|_{L^2([0,\ 1],\ H^1(\R^2))}^2\leq C\,\eps^{2-2\alpha-\beta},\quad \|u_{a,\eps}\|_{L^\infty([0,\ 1]\times\R^2)}\leq \eps^{1-\alpha},
\end{equation}
which implies for the flow associated
\[
\partial_t X(t,0,x)=u_{a,\eps}(t,X(t,0,x)),
\]
that there exists a set $\Omega_{a,\eps}\subset B(a,\eps)$ with $|\Omega_{a,\eps}|\geq \frac{1}{4}\,|B(a,\eps)|$, $c_\gamma>0$ and $0<\delta_\eps<\eps$ s.t. if $x,\;y\in \Omega_{a,\eps}$ with $|x-y|\leq \delta_\eps$ then for any $t\leq 1$
\begin{equation}
\E(|X(t,0,x)-X(t,0,y)|^\gamma)\geq \frac{(1+c_\gamma)^{t\,\eps^{-\alpha}/C}}{2}\,|x-y|.\label{expectation} 
\end{equation}\label{randomprop}
\end{prop}
Note that the flow is random only because $u_{a,\eps}$ is.

Prop. \ref{randomprop} is a consequence of a stretching flow combined with appropriate mixing in orientation. The stretching part is mostly a 
rotation by varying angle. More precisely it belongs to the following class
\begin{defi}
Denote by ${\cal R}_a$ the class of non linear transform $F:\R^2\rightarrow \R^2$ s.t. 
\[
\forall r>0,\exists \theta(r),\quad \mbox{if}\ |x-a|=r\qquad F(x)=a+e^{i\,\theta(r)}\,(x-a),
\]
where the classical complex notation is used for rotations.
\end{defi}
The stretching structure itself is given by
\begin{lemma}
For any $\tau$ small enough, $a\in\R^2$, any $\gamma>0$, $\eps>0$, $\alpha<1$, there exists a time-independent, divergence free vector field $u_{a,\eps}^r$ with support in the ball $B(a,\;\eps)$ and
\begin{equation}
\|u_{a,\eps}^r\|_{H^1(\R^2)}^2\leq C\,\eps^{2-2\alpha},\quad \|u_{a,\eps}^r\|_{L^\infty([0,\ 1]\times\R^2)}\leq \eps^{1-\alpha},\label{H1ur}
\end{equation}
which implies that the flow associated 
\[
\partial_t X(t,t_0,x)=u_{a,\eps}^r(X(t,t_0,x)),\quad X(t=t_0,t_0,x)=x,
\]
belongs to ${\cal R}_a$, and that there exists $C>0$ and $l<1$, for any $x,\;y\in B(a,\eps)\setminus B(a,\eps/2)$ 
\begin{equation}
\E(|X(t_0+\eps^\alpha,t_0,x)-X(t_0+\eps^\alpha,t_0,x)|^\gamma)\geq |x-y|^\gamma\,\left(1+C\,\rho^4\,\omega_1^2\,(1-2\,l\,\omega_2^2\right),
\label{expstretch}
\end{equation}
where
\[
\rho=\frac{|x-a|}{\eps},\quad \omega_1=\frac{(x-a)\cdot(x-y)}{|x-a|\,|x-y|},\quad \omega_2=\frac{(x-a)^\perp\cdot(x-y)}{|x-a|\,|x-y|}.
\]
In addition almost surely for every $x$ and $y$
\[
|X(t_0+\eps^\alpha,t_0,x)-X(t_0+\eps^\alpha,t_0,x)|\leq C\,|x-y|,
\]
and if $|x-a|\leq \eps/2$ then $X(t_0+\eps^\alpha,t_0,x)$ is just the image from $x$ by a rotation of angle $\tau$.  
\label{lemmastretching}
\end{lemma} 
If $\omega=(\omega_1,\;\omega_2)$ could automatically be chosen randomly and uniformly on the unit circle then the expectation of $\omega_1^1\,(1-2l\,\omega_2^2)$ would be strictly positive as
\[
\oint_{S^1} \omega_1^1\,(1-2\,l\,\omega_2^2)\,d\omega=\frac{1}{2}-\frac{l}{2}>0.
\]
Thus the second step in the construction is to add enough rotations to ensure that the orientation is close to being uniformly distributed per
\begin{lemma}
For any $\tau$ s.t. $2\pi/\tau\in\N$, any $a\in\R^2$, any $\eps>0$, $0<\alpha<1$, $\beta>\alpha$, and $\eta>0$, there exists $C>0$ and a random divergence free vector field $R_{a,\eps}^{t_0}$ with support in the ball $B(a,\;\eps)$ and
\begin{equation}
\|R_{a,\eps}^{t_0}\|_{L^2([t_0,\ t_0+2\eps^{\alpha}],\ H^1(\R^2))}^2\leq C\,\eps^{\alpha}\,\eps^{2-2\alpha-\beta}, \quad \|R_{a,\eps}^{t_0}\|_{L^\infty([0,\ 1]\times\R^2)}\leq \eps^{1-\alpha}.\label{H1Rt}
\end{equation}
That field implies for the flow associated
\[
\partial_t X(t,t_0,x)=R^{t_0}_{a,\eps}(t,X(t,t_0,x)),\quad X(t=t_0,t_0,x)=x,
\]
that there exists a set $A_{a,\eps}$ with $|A_{a,\eps}|\leq C\,\eps^{1+\beta}$
s.t. 
\begin{itemize}
\item[i.] Almost surely
\[
\frac{|x-a|}{\eps}-\eta\leq \frac{|X(t_0+2\,\eps^\alpha,t_0,x)-a|}{\eps}\leq \frac{|x-a|}{\eps}+\eta.
\]
\item[ii.] Define as $O_{a,\eps}$ the smallest invariant set by the flow containing the annulus $\{\eps/2\leq|x-a|\leq \eps\}$. Then $O_{a,\eps}$ is included in $B(a,2\,\eps)$ and moreover $B(a,\eps/2)\cap O_{a,\eps}$ is invariant by rotation centered at $a$ and of angle $\tau$.  
\item[iii.] If $x,\;y\in O_{a,\eps}\setminus A_{a,\eps}$ and $|x-y|\leq \eps^{1+\beta}$, both $X(t_0+2\,\eps^\alpha,t_0,x)$ and $X(t_0+2\,\eps^\alpha,t_0,y)$ belong to $\{\eps/2\leq|x-a|\leq \eps\}$ with probability at least $1/2$. 
\item[iv.] For any fixed $F\in {\cal R}_a$, the probability that $X(t_0+2\,\eps^\alpha,t_0,x)$ or $X(t_0+2\,\eps^\alpha,t_0,y)$ belongs to $F\,A_{a,\eps}$ is at most $C\,\eps^\beta$.
\item[v.] If $x,\;y\in O_{a,\eps}\setminus A_{a,\eps}$, $|x-y|\leq \eps^{1+\beta}$ then 
\[
|X(t_0+2\,\eps^\alpha,t_0,x)-X(t_0+2\,\eps^\alpha,t_0,y)|=|x-y|,
\]
with probability at least $1-C\,\eps^\beta$.
\item[vi.] If $x,\;y\in O_{a,\eps}\setminus A_{a,\eps}$ and $|x-y|\leq\eps^{1+\beta}$, then with probability at least $1-2\,\eta$, 
\[
\frac{X(t_0+2\,\eps^\alpha,t_0,x)-X(t_0+2\,\eps^\alpha,t_0,y)}{|X(t_0+2\,\eps^\alpha,t_0,x)-X(t_0+2\,\eps^\alpha,t_0,y)|}
\]
is {\bf uniformly} randomly distributed on $S^1$. 
\end{itemize}
\label{lemmamixing}
\end{lemma}
Remark that contrary to the field $u_{a,\eps}^r$ in Lemma \ref{lemmastretching}, the field $R_{a,\eps}^{t_0}$ is time-dependent. As seen in the construction it is actually the combination of two time-independent fields.   
\subsection{From Lemmas  \ref{lemmastretching} and \ref{lemmamixing} to Prop. \ref{randomprop}}
We combine the fields obtained from Lemmas \ref{lemmastretching} and \ref{lemmamixing}. Denote $t_k=k\,\eps^\alpha$, $N$ s.t. $t_{3N}\sim 1$ (see the end of the proof for a more precise bound). Take $N$ independent copies $u_{a,\eps}^{r,n}$ and $R_{a,\eps}^{t_{3n}}$, $n=0\ldots N-1$ of the fields from Lemmas \ref{lemmastretching} and \ref{lemmamixing} with the same parameters $a$, $\eps$, $\alpha$, $\beta$ and $\tau$ s.t. $2\,\pi/\tau\in\N$. The parameter $\eta$ will be chosen later.

Define 
\begin{equation}
u_{a,\eps}=\sum_{n=0}^{N-1} \Big\{ \ind_{t\in[t_{3n},\ t_{3n+2})}\,R_{a,\eps}^{t_{3n}} + \ind_{t\in[t_{3n+2},\ t_{3n+3})}\,u_{a,\eps}^{r,n}\Big\}.\label{defuaeps}
\end{equation}
%
By the properties of all the fields in the sum, the field $u_{a,\eps}$ is divergence free, supported in $B(a,2\,\eps)$. In addition combining the estimates \eqref{H1ur} and \eqref{H1Rt} over the time interval $[0,\ 1]$, one checks that
\[
\|u_{a,\eps} \|^2_{L^2([0,\ 1],\ H^1(\R^2))}\leq C\,\eps^{2-2\alpha-\beta}.
\]
Note by the way that the contribution from the $R_{a,\eps}^{t_{3n}}$ is the dominant one.

As before, denote by $X$ the flow from $t=0$ to $t_{3N}$
\[
\partial_t X(t,s,x)=u_{a,\eps}(t,X(t,0,x)),\quad X(t,s,x)=x.
\] 
The sequence obtained by taking the flow, $X(t_{3n},0,x)$ and $X(t_{3n+2},0,x)$, at the times $t_{3n}$ and $t_{3n+2}$  constitutes a {\em discrete Markov chain}, simply by the independence of all the $u_{a,\eps}^{r,n}$ and $R_{a,\eps}^{t_{3n}}$. 

First we prove by induction that if $x\in O_{a,\eps}$ then  $X(t_{3n},0,x)\in O_{a,\eps}$ almost surely for any $n$. 

Assume that $X(t_{3n},0,x)\in O_{a,\eps}$, one wishes to show that $X(t_{3n+3},0,x)\in O_{a,\eps}$.
Of course by its definition $O_{a,\eps}$ is invariant by the flow of any $R_{a,\eps}^{t_{3n}}$ and thus $X(t_{3n+2},0,x)\in O_{a,\eps}$. By its definition $\{\eps/2\leq |x-a|\leq \eps\}\subset O_{a,\eps}$ and the annulus is invariant by the flow of $u_{a,\eps}^{r,n}$. Since $u_{a,\eps}^{r,n}$ is compactly supported in $B(a,\eps)$, it also leaves invariant $O_{a,\eps}\setminus B(a,\eps)$. Finally if  $X(t_{3n+2},0,x)\in O_{a,\eps}\cap B(a,\eps/2)$ then by Lemma \ref{lemmastretching}, as $t_{3n+3}=t_{3n+2}+\eps^\alpha$, $X(t_{3n+3},0,x)$ is the image of $X(t_{3n+2},0,x)$ by a rotation of angle $\tau$. By Lemma \ref{lemmamixing}, $O_{a,\eps}\cap B(a,\eps/2)$ is invariant by such rotation. As a consequence $X(t_{3n+3},0,x)\in O_{a,\eps}$.

We now estimate, again by induction, the total probability that a trajectory $X(t,0,x)$ would enter the bad set $A_{a,\eps}$ at any one time $t=t_{3n}$. The key is the fact that the $X(t_{3n},0,x)$ are a Markov chain through the independence of each field in the definition \eqref{defuaeps} of $u_{a,\eps}$. 

Assume that $X(t_{3n},0,x)\not\in A_{a,\eps}$. For any realization of $u_{a,\eps}^{r,n}$, the flow $X(t_{3n+3},t_{3n+2},.)$ and its inverse $X(t_{3n+2},t_{3n+3},.)$ are transforms of ${\cal R}_a$, fixed in the sense that they are independent of the realization of $R_{a,\eps}^{t_{3n}}$. Therefore by Lemma \ref{lemmamixing}, $X(t_{3n+2},0,x)\not\in X(t_{3n+2},t_{3n+3}, A_{a,\eps})$ with probability at least $1-C\,\eps^{\beta}$. By the independence of the fields, one hence has that
\[
\mathbb{P} \left(X(t_{3n+3},0,x)\in A_{a,\eps}\;|\ X(t_{3n},0,x)\not\in A_{a,\eps}\right)\leq C\,\eps^{\beta}.
\]
Note that the same is true if one considers two trajectories starting from $x$ and $y$ at the same time by Lemma \ref{lemmamixing}. By the Markov chain property, if $x,\;y\not\in A_{a,\eps}$ then 
\begin{equation}
\mathbb{P} \left(X(t_{3n},0,x)\not\in A_{a,\eps}\ \mbox{and}\ X(t_{3n},0,y)\not\in A_{a,\eps},\quad\forall n\right)\geq 1-C\,N\,\eps^{\beta}.\label{badsetproba}
\end{equation}
Still considering two trajectories $X(t,0,x)$ and $X(t,0,y)$, denote for convenience 
\[
\delta_k=X(t_k,0,x)-X(t_k,0,y).
\]
Take $x$ and $y$ close enough initially so that $|\delta_k|< \eps^{1+\beta}$ for any $k\leq 3\,N$ with large probability. By the previous point, with probability $1-C\,N\,\eps^{\beta}$, we can assume that $X(t_{3n},0,x)\not\in A_{a,\eps}$ and $X(t_{3n},0,y)\not\in A_{a,\eps}$. Thus assuming that $|\delta_{3n}|<\eps^{1+\beta}$, one has that $|\delta_{3n+2}|=|\delta_{3n}|$ by Lemma \ref{lemmamixing}. Then by Lemma \ref{lemmastretching}, $|\delta_{3n+3}|\leq C\,|x-y|$.

This implies that if
\begin{equation}
|x-y|< \eps^{1+\beta}\,C^{-N},\label{x-y}
\end{equation} 
then $|\delta_k|<\eps^{1+\beta}$ for all $k\leq 3\,N$ with probability at least $1-C\,N\,\eps^{\beta}$. Therefore we choose for Prop. \ref{randomprop}
\[
\Omega_{a,\eps}=B(a,\eps)\cap O_{a,\eps}\setminus A_{a,\eps}, \quad \delta_\eps=\eps^{1+\beta}\,C^{-N}.
\]
Denote by ${\cal F}_k$ the filtration adapted to the random fields until time $t_k$. 
Finally denote by $\E_C$ the expectation conditioned to having both $X(t_{3n},0,x)\not\in A_{a,\eps}$ and $X(t_{3n},0,y)\not\in A_{a,\eps}$ for every $n\leq n$. Write
\[\begin{split}
&\rho_k=\frac{|X(t_k,0,x)-a|}{\eps},\quad\omega_k=(\omega_{1,k},\;\omega_{2,k}),\\
&\omega_{1,k}=\frac{\delta_k}{|\delta_k|}\cdot \frac{X(t_k,0,x)-a}{|X(t_k,0,x)-a|},\quad 
 \omega_{2,k}=\frac{\delta_k}{|\delta_k|}\cdot \frac{(X(t_k,0,x)-a)^\perp}{|X(t_k,0,x)-a|}.
\end{split}\]
By \eqref{expstretch} of Lemma \ref{lemmastretching}, if both $X(t_{3n+2},0,x)\in \{\eps/2\leq |x-a|\leq \eps\}$ and $X(t_{3n+2},0,y)\in \{\eps/2\leq |x-a|\leq \eps\}$
\[\begin{split}
&\mathbb{E}\left(|\delta_{3n+3}|^\gamma\;|\ {\cal F}_{3n+2}\right)\geq |\delta_{3n+2}|^\gamma\,\left(1+2\gamma\,\tau^2\,\rho_{3n+2}^4 \,\left|\omega_{1,3n+2}\right|^2\,(1-2l\,\left|\omega_{2,3n+2}\right|^2)\right).\\
\end{split}
\]
On the other hand if $X(t_{3n+2},0,x)\not\in \{\eps/2\leq |x-a|\leq \eps\}$ or 
$X(t_{3n+2},0,y)\not\in \{\eps/2\leq |x-a|\leq \eps\}$ then
\[
|\delta_{3n+3}|^\gamma=|\delta_{3n+2}|^\gamma.
\]
Assume now that $X(t_{3n},0,x)\in O_{a,\eps}\setminus A_{a,\eps}$ and $X(t_{3n},0,y)\in O_{a,\eps}\setminus A_{a,\eps}$. Apply Lemma \ref{lemmamixing} which in particular implies that with probability at least $1/2$, both $X(t_{3n+2},0,x)\in \{\eps/2\leq |x-a|\leq \eps\}$ and $X(t_{3n+2},0,y)\in \{\eps/2\leq |x-a|\leq \eps\}$, thus
\[\begin{split}
&\mathbb{E}\left(|\delta_{3n+3}|^\gamma\;|\ {\cal F}_{3n}\right)\\
&\qquad\geq \mathbb{E}\left(|\delta_{3n+2}|^\gamma\,(1+\gamma\,\tau^2\,\rho_{3n+2}^4 \,\left|\omega_{1,3n+2}\right|^2\,(1-2l\,\left|\omega_{2,3n+2}\right|^2))\;|\ {\cal F}_{3n}\right).\\
\end{split}
\]
Still by Lemma \ref{lemmamixing}, with probability at least $1-C\eps^\beta$, one has that $|\delta_{3n+2}|=|\delta_{3n}|$ which gives
\[\begin{split}
&\mathbb{E}\left(|\delta_{3n+3}|^\gamma\;|\ {\cal F}_{3n}\right)\\
&\quad\geq |\delta_{3n}|^\gamma\,(1-C\,\eps^\beta)\,\mathbb{E}\left(1+\gamma\,\tau^2\,\rho_{3n+2}^4 \,\left|\omega_{1,3n+2}\right|^2\,(1-2l\,\left|\omega_{2,3n+2}\right|^2)\;|\ {\cal F}_{3n}\right).\\
\end{split}
\]
Apply a last time Lemma \ref{lemmamixing}, point $i.$ to deduce that $\rho_{3n+2}=\rho_{3n}\pm\eta\,\eps$ and point $vi.$ to deduce that with probability at least $1-\eta$, $\delta_{3n+2}/|\delta_{3n+2}|$ is uniformly distributed on $S^1$ and hence by definition so is $\omega_{3n+2}$. Thus
\[\begin{split}
&\mathbb{E}\left(1+\gamma\,\tau^2\,\rho_{3n+2}^4 \,\left|\omega_{1,3n+2}\right|^2\,(1-2l\,\left|\omega_{2,3n+2}\right|^2)\;|\ {\cal F}_{3n}\right)\\
&\quad\geq 1-C\,\eps^4\,\eta+\,\gamma\,\tau^2\,\rho_{3n}^4\,\int_{S^1} |\omega_1|^1\,(1-2l\,|\omega_2|^2)\,d\omega.
\end{split}
\]
Since $l<1$ then $\int_{S^1} |\omega_1|^1\,(1-2l\,|\omega_2|^2)\,d\omega>0$. Choose now $\eta$ small enough with respect to $\gamma$, $\tau$ and $1-l$ to find that for some $c_\gamma>0$
\[\begin{split}
&\mathbb{E}\left(1+\gamma\,\tau^2\,\rho_{3n+2}^4 \,\left|\omega_{1,3n+2}\right|^2\,(1-2l\,\left|\omega_{2,3n+2}\right|^2)\;|\ {\cal F}_{3n}\right)\geq 1+c_\gamma.
\end{split}
\]
This leads to the fact that if $X(t_{3n},0,x)\in O_{a,\eps}\setminus A_{a,\eps}$ and $X(t_{3n},0,y)\in O_{a,\eps}\setminus A_{a,\eps}$ then
\begin{equation}
\mathbb{E}\left(|\delta_{3n+3}|^\gamma\;|\ {\cal F}_{3n}\right)\geq (1+c_\gamma)\,|\delta_{3n}|^\gamma.\label{basisinduction}
\end{equation}
Repeating by induction \eqref{basisinduction}, one finds that
\[
\mathbb{E}_C\left(|\delta_{3n}|^\gamma\right)\geq |x-y|\,(1+c_\gamma)^n.
\]
Finally since all $X(t_{3n},0,x)\in O_{a,\eps}\setminus A_{a,\eps}$ and $X(t_{3n},0,y)\in O_{a,\eps}\setminus A_{a,\eps}$ with probability $1-C\,N\,\eps^\beta$ by \eqref{badsetproba}
\[
\mathbb{E}\left(|\delta_{3n}|^\gamma\right)\geq (1-C\,N\,\eps^\beta)\,|x-y|\,(1+c_\gamma)^n.
\] 
It only remains to choose $N$ s.t. $C\,N\,\eps^\alpha\leq 1/2$ (recall that $\beta>\alpha$) and note that $n=t_{3n}/(3\,\eps^\alpha)$ to conclude the proof.
\subsection{From Prop. \ref{randomprop} to Theorem \ref{counterexample}}
The first step is to remove the randomness by finding one realization of the random field $u_{a,\eps}$ which produces a growth like \eqref{expectation} for enough $x$ and $y$. 

Apply Prop. \ref{randomprop} and average \eqref{expectation} for every $x\in B(a,\eps)$ and every $y$ in $B(x,\,\delta_\eps)$ 
\[\begin{split}
&\oint_{B(a,\eps)}\oint_{B(x,\,\delta_\eps)}\mathbb{E}\left(|X(1,0,x)-X(1,0,y)|^\gamma\right)\,dy\,dx\\
&\qquad\geq \frac{1}{4}\oint_{\Omega_{a,\eps}}\oint_{B(x,\,\delta_\eps)}\mathbb{E}\left(|X(t_{3N},0,x)-X(t_{3N},0,y)|^\gamma\right)\,dy\,dx\\
&\qquad\qquad\geq \frac{(1+c_\gamma)^N}{C}\,\delta_\eps^\gamma,
\end{split}\]
for some $C>0$, as $|\Omega_{a,\eps}|\geq \frac{1}{4}|B(a,\eps)|$ and $t_{3N}\leq 1$.

Therefore there exists at least one realization ({\em i.e.} deterministic), denoted by $u_{a,\eps}^d$, of the random field $u_{a,\eps}$ s.t. 
when $X$ solves
\[
\partial_t X=u^d_{a,\eps}(t,X),\quad X(t=0,x)=x,
\]
then
\begin{equation}\begin{split}
&\oint_{B(a,\eps)}\oint_{B(x,\,\delta_\eps)}\left(|X(1,0,x)-X(1,0,y)|^\gamma\right)\,dy\,dx\\
&\qquad\qquad\geq \frac{(1+c_\gamma)^{C\,\eps^{-\alpha}}}{C}\,\delta_\eps^\gamma. \label{detergrowth}
\end{split}\end{equation}
Note that there are of course many such realizations but one is enough for our purpose.

The next step is to combine several such structures. Choose $\alpha$ and $\beta>\alpha$ s.t. $2-2\,\alpha-\beta>0$; remark that indeed $N\,\eps^\beta\sim \eps^{\beta-\alpha}<<1$. Choose centers $a_i$, $i\in\N$, s.t. $\{a_i,\ i\in\N\}$ is compactly supported and
\[
B(a_i,\,2^{-i+1})\cap B(a_j,\,2^{-j+1})=\emptyset,\qquad \mbox{if}\ i\neq j,
\] 
which is always possible as the total volume of the balls $B(a_i,\,2^{-i+1})$ is finite. Finally define
\begin{equation}
u(t,x)=\sum_{i=0}^\infty u^d_{a_i,2^{-i}}.\label{defu}
\end{equation}
The $u^d_{a_l,2^{-l}}$ have compact support with empty intersections. Therefore $u$ is divergence free, $\|u\|_{L^\infty}\leq 1$ since $\alpha<1$ and
\[
\|u\|_{L^2([0,\ 1],\ H^1(\R^2))}^2=\sum_{i=0}^{+\infty} \|u^d_{a_i,2^{-i}}\|_{L^2([0,\ 1],\ H^1(\R^2))}^2\leq C\,\sum_{i=0}^{+\infty} \left(2^{-i}\right)^{2-2\alpha-\beta}<\infty.
\]
Let $X$ solve for any $x$
\[
\partial_t X(t,0,x)=u(t,X(t,0,x)),\quad X(0,0,x)=x.
\]
The dynamics is resolved independently in any $B(a_i,\;2^{-i+1})$ and in particular it is well posed as $u$ is Lipschitz in any such ball (but of course with a Lipschitz constant depending on $i$). However $X(t,0,x)$ does not belong to any Sobolev space $W^{1,\gamma}$ because
\[
\limsup_{h\rightarrow 0} \int_{\R^2} \oint_{B(x,h)}\frac{|X(1,0,x)-X(1,0,y)|^\gamma}{h^\gamma}\,dy\,dx=+\infty.
\]
Indeed if $\delta_{2^{-i}}\leq h<2\,\delta_{2^{-i}}$ then
\[\begin{split}
&\int_{\R^2} \oint_{B(x,h)}\frac{|X(1,0,x)-X(1,0,y)|^\gamma}{h^\gamma}\,dy\,dx\\
&\qquad\qquad\geq
C\,\int_{B(a,2^{-l})} \oint_{B(x,\delta_{2^{-l}})}\frac{|X(t_{3N},0,x)-X(t_{3N},0,y)|^\gamma}{h^\gamma}\,dy\,dx\\
&\qquad\qquad\geq C\,2^{-2l}\,(1+c_\gamma)^{C\,2^l},
\end{split}\]
by \eqref{detergrowth}. This last quantity tends to $+\infty$ as $i$ tends to $+\infty$ and hence $h$ to $0$, which concludes the proof of Theorem \ref{counterexample}.

\bigskip

$\bullet$ {\bf Final remarks.} The choice of $u^d_{a,\eps}$ and then $u$ may be deterministic but it relies on the choices of all the parameters entering in the construction of $u_{a,\eps}^r$ and $R_{a,\eps}^{t_0}$. Those choices are not explicit and very likely not obvious at all. This is the main advantage of using randomness.

There are many possible choices of $\alpha$ and $\beta$, for instance $\beta=1/2$, $\alpha=1/4$. Ideally one would like to take $\beta=+\infty$ and $\alpha=1$ however $u_{a,\eps}$ would not belong to $H^1_x$ but only to $BV$. 
While $\beta=+\infty$ would seem to be a problem as it enters for instance in the definition of $\delta_\eps$, the analysis could carried out with the less demanding
\begin{equation}
|x-y|<< \eps\,e^{-C\,\eps^{-\alpha}}.\label{x-yrelax}
\end{equation}
\section{The stretching part: Proof of Lemma \ref{lemmastretching}}
\subsection{The basic deterministic stretching} 
For any $\eps>0$ and any point $a\in \R^2$, define
\[
u_{a,\eps}^s= \eps^{-\alpha}\,f_\eps(|x-a|^2)\,(x-a)^{\perp},
\]
where $(x-a)^\perp$ denotes the rotation of the vector $x-a$ by $\pi/2$, and $f_\eps$ is given by
\[
f_\eps(\xi)=\tau'\,(1-\xi/\eps^2)\,\ind_{\eps^2/4< \xi< \eps^2}+\tau'\,(1-1/4)\,\ind_{\xi\leq \eps^2/4}.
\]
%
Note that $u_{a,\eps}^s$ is Lipschitz with support in the ball $B(a,\eps)$. Moreover it is divergence free
\[
\div u^s_{a,\eps}=\eps^{-\alpha}\,f_\eps'(|x-a|^2) \,(-(x_1-a_1)\,(x_2-a_2)+(x_2-a_2)\,(x_1-a_1))=0.
\]
Finally 
\begin{equation}
\|u_{a,\eps}^s\|_{H^1_x}\sim \eps^{2-2\alpha}, \|u_{a,\eps}^s\|_{L^\infty}\leq \eps^{1-\alpha}.\label{H1us}
\end{equation}
The solution from $t=t_0$ till $t=t_1=t_0+\eps^\alpha$ to 
\[
\partial_t X(t,s,x)=u_{a,\eps}^s(t,X(t,s,x)),\quad X(t=s,s,x)=x
\]
is simple enough: If $|X(t_0,s,x)-a|>\eps$ then $X(t,s,x)=X(t_0,s,x)$. If $|X(t_0,s,x)-a|\leq\eps/2$ then the trajectory is not constant but at $t=t_1$ it will have rotated by exactly $3\tau'/4$. 

In the general case, $|X-a|$ is constant in time and 
\[
X(t_1,s,x)=a+e^{i\,f_\eps(|X(t_0,s,x)-a|^2)}\,(X(t_0,s,x)-a),
\]
where the complex exponential is used to denote the corresponding rotation in $\R^2$.

Now compare two trajectories $X(t,s,x)$ and $X(t,s,y)$. Denote
\[
\delta=X(t_0,s,y)-X(t_0,s,x).
\] 
Assuming that initially $X(t_0,s,y)$ and $X(t_0,s,x)$ are in the annulus $\{\eps/2< |x-a|< \eps\}$, then 
\begin{equation}
\begin{split}
&|X(t_1,s,x)-X(t_1,s,y)|^2\\
&\qquad =\left|e^{i\,\tau'\,\eps^{-2}\,\delta\cdot (X(t_0,s,x)+X(t_0,s,y)-2a)}\,(X(t_0,s,x)-a)-Y(t_0,s,x)+a  \right|^2\\
&\qquad=\left|\delta+2\,\tau'\,\eps^{-2}\,\delta\cdot (X(t_0,s,x)-a)\,(X(t_0,s,x)-a)^\perp\right|^2+O(\delta^3\,\eps^{-1}).
\end{split}\label{deltaX}
\end{equation}
Depending on the angle between $\delta$ and $X(t_0,s,x)$, the distance may be multiplied by up to $1+2\,\tau'$ but it can also decrease by a similar amount. 
\subsection{Random sense of rotation}
%
Randomize the previous field by picking a random sense of rotation leading to
\[
u_{a,\eps}^r=\sigma\,u_{a,\eps}^s(x),
\]
where $\sigma$ is $+1$ with probability $1/2$ and $-1$ otherwise.
%
By \eqref{H1us}, the field $u_{a,\eps}^r$ satisfies \eqref{H1ur}.

From $t=t_0$ till $t=t_1=t_0+\eps^\alpha$ let $X$ be the flow and solve 
\[
\partial_t X(t,t_0,x)=u_{a,\eps}^r(t,X(t,t_0,x)), \quad X(t=t_0,t_0,x)=x.
\]
Compare again two trajectories $X(t,s,x)$ and $X(t,s,y)$ for $x$ and $y$ very close. Denote still $\delta=X(t_0,s,y)-X(t_0,s,x)=x-y$ ($\delta$ is deterministic here). To finish the proof of Lemma \ref{lemmastretching}, one has to get a lower bound for $\mathbb{E} |X(t_1,s,x)-X(t_1,s,y)|^\gamma$ for any $\gamma>0$.

Denote as in the Lemma
\[\begin{split}
&\rho=\frac{|X(t_0,s,x)-a|}{\eps},\quad \omega_1=\frac{\delta}{|\delta|}\cdot \frac{X(t_0,s,x)-a}{|X(t_0,s,x)-a|},\\
& \omega_2=\frac{\delta}{|\delta|}\cdot \frac{(X(t_0,s,x)-a)^\perp}{|X(t_0,s,x)-a|}.
\end{split}\]
Following \eqref{deltaX}, calculate for $\eps/2\leq |x-a|\leq \eps$ and $\eps/2\leq |y-a|\leq \eps$, that is $1/2\leq\rho\leq 1$
\[\begin{split}
&\sum_{\sigma=\pm 1}\big|\delta
+2\,\tau'\,\sigma\eps^{-2}\,\delta\cdot (X(t_0,s,x)-a)\,(X(t_0,s,x)-a)^\perp\big|^\gamma\\
&\ =\sum_{\sigma=\pm 1}\,|\delta|^\gamma\left(1+4\tau'^2\,\rho^4\,\omega_1^2 +4\sigma\,\tau'\,\rho^2\,\omega_1\,\omega_2
\right)^{\gamma/2}.
\end{split}\]
Note that there exists a universal constant $C$ s.t. for any $0<\gamma<1$ and any $-1<\xi<1$
\[
(1+\xi)^{\gamma/2}\geq 1+\frac{\gamma}{2}\,\xi-\frac{\gamma\,(1-\gamma/2)}{4}\,\xi^2 -C\,\gamma\,|\xi|^3.
\]
Thus
\[\begin{split}
&|\delta|^{-\gamma}\,\sum_{\sigma=\pm 1}\big|\delta
+2\,\tau'\,\sigma\eps^{-2}\,\delta\cdot (X(t_0,s,x)-a)\,(X(t_0,s,x)-a)^\perp\big|^\gamma\\
&\ \geq 2+\frac{\gamma}{2}\sum_{\sigma=\pm1} (4\tau'^2\,\rho^4\,\omega_1^2 +4\sigma\,\tau'\,\rho^2\,\omega_1\,\omega_2-8\,(1-\gamma/2)\tau'^2\,\rho^4\,\omega_1^2\,\omega_2^2)\\
&\qquad+O(\gamma\,\tau'^3\,\rho^6\,\omega_1^3).\\
\end{split}
\]
Now for $\gamma>0$ fixed, by taking $\tau'$ small enough, there exists $l<1$ s.t.
\begin{equation}\begin{split}
&|\delta|^{-\gamma}\,\sum_{\sigma=\pm 1}\big|\delta
+2\,\tau'\,\sigma\eps^{-2}\,\delta\cdot (X(t_0,s,x)-a^k)\,(X(t_0,s,x)-a^k)^\perp\big|^\gamma\\
&\ \geq 2+4\,\gamma\,\tau'^2\,\rho^4\,\omega_1^2\,(1-2\,l\,\omega_2^2).
\end{split}\label{deltaomega}
\end{equation}
This concludes the proof of Lemma \ref{lemmastretching}, taking $\tau=3\,\tau'/4$ or $\tau'=4\,\tau/3$.
%
%
\section{The orientation mixing structure: Proof of Lemma \ref{lemmamixing}} 
Define the field
\[
R_{z,\eps,r}=\eps^{-\alpha}\, g_\eps((x-z)/(r\eps))\,(x-z)^\perp,
\]  
where $r$ will be chosen later in terms of $\eta$. The function $g_\eps$ is radially symmetric, smooth with compact support in $B(0,1)$ and such that $g_\eps(x)=1$ in the ball $ |x|\leq 1-\eps^{\beta}$.

We also define the corresponding field in the annulus
\[
R_{z,\eps,r}^{ann}=\eps^{-\alpha}\, g_\eps^{ann}((x-z)/(r\eps))\,(x-z)^\perp,
\]
where  $g_\eps^{ann}$ is radially symmetric, smooth with compact support in $\{1/2\leq |x|\leq 1\}$ and such that $g_\eps^{ann}(x)=1$ in the annulus $1/2+\eps^\beta\leq |x|\leq 1-\eps^{\beta}$.

Therefore $R_{z,\eps,r}$ generates a rotation in  $B(z,r\,\eps(1-\eps^\beta))$, has compact support in $B(z,r\,\eps)$, and is a more complicated transform only in the thin annulus 
\[
A_{z,\eps,r}=\{r\,\eps(1-\eps^\beta)<|z-x|<r\,\eps\}.
\]
Since $g_\eps$ is radially symmetric then $R_{z,\eps,r}$ is divergence free. Moreover since $g_\eps$ is smooth, $|\nabla g_\eps|\sim \eps^{-\beta}$ on the annulus $A_{z,\eps,r}$  and vanishes elsewhere, we have 
\begin{equation}
\|R_{z,\eps,r}\|_{L^2([T_1,\ T_2],\;H^1_x)}^2\sim (T_2-T_1)\,\eps^{2-2\alpha-\,\beta},\quad \|R_{z,\eps,r}\|_{L^\infty([0,\ 1]\times\R^2)}\leq \eps^{1-\alpha}.\label{H1Rz}
\end{equation}
Now for any $\eta>0$, one can find a finite number of centers $z_i$, $i\leq k$, and radii $r_i<\eta$ s.t. the number of spheres, radii $r_i$ and $(z_i-a)/\eps$ are independent of $\eps$ and
\begin{itemize}
\item[i.] The spheres are disjoint: $B(z_i,r_i\,\eps)\cap B(z_j,r_i\,\eps)=\emptyset$, $\forall i\neq j$.
\item[ii.] $\bigcup_{i=1...k} B(z_i,r_i\,\eps)\subset \{(1/2-\eta)\,\eps<|x-a|<(1+\eta)\,\eps\})$ but
\[
|B(z_i,r_i\,\eps)\cap \{\eps/2<|x-a|<\eps\}|\geq \frac{1}{4}\,|B(z_i,r_i\,\eps)|.
\]
\item[iii.] For any $\eps/2<R<\eps$, most of the circle $\{|x-a|=R\}$ is included in $\bigcup_{i=1...k} B(z_i,r_i\,\eps\,(1-\eps^\beta))$. More precisely
\[
\mbox{length}\left(\{|x-a|=R\}\setminus \bigcup_{i=1...k} B(z_i,r_i\,\eps\,(1-\eps^\beta))\right)<\eta\,R.
\] 
\item[iv.] For convenience we assume that the balls $B(z_i,r_i\,\eps)$ intersecting with $B(a,\eps/2)$ are all the image of one another by rotations of angle $\tau$; this is possible as  $2\,\pi/\tau\in \N$.
\end{itemize}
Note that because of $iii$, it would not possible to have $\bigcup_{i=1...k} B(z_i,r_i\,\eps)\subset \{\eps/2<|x-a|<\eps\}$, hence the more convoluted version of $ii$.
\begin{figure}[!htb]
\centering
\includegraphics[scale=.5]{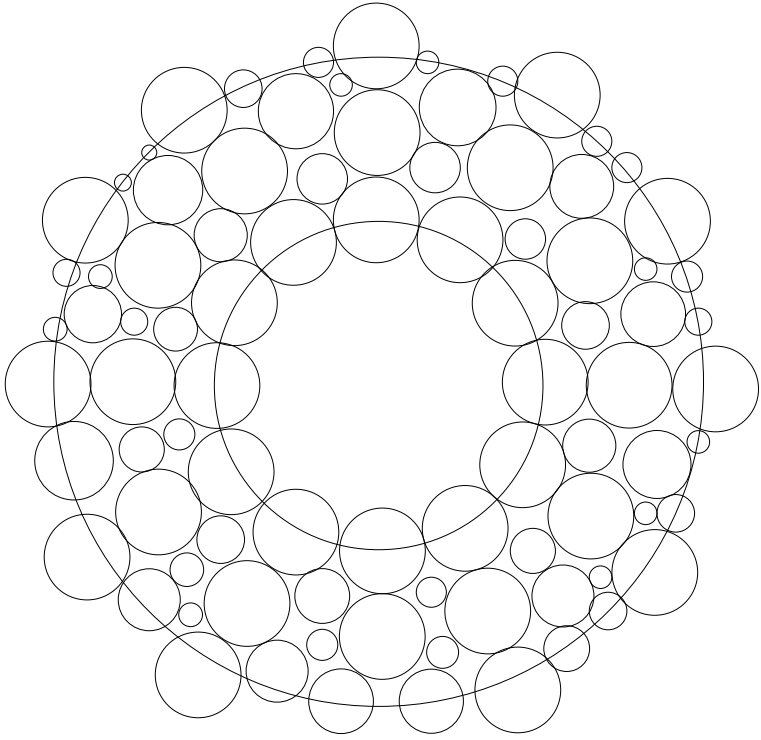}
\caption{An example of choice of spheres}
\label{covering}
\end{figure}

\bigskip

Starting from any time $t_0$, one achieves {\em random mixing} by combining many smaller rotations by defining
\[
R_{a,\eps}^{t_0}=\theta\,R_{a,\eps,1}^{ann}\,\ind_{t\in[t_0,\ t_0+\eps^\alpha]}+\ind_{t\in(t_0+\eps^\alpha,\ t_0+2\,\eps^\alpha]}\sum_{i=1}^k \lambda_i R_{z_i,\eps,r_i},
\] 
where each $\lambda_i$ is an independent random number uniformly chosen in $[0,\ 2\pi]$ and so is the random angle $\theta\in [0,\ 2\pi]$.

By \eqref{H1Rz} and since the fields have disjoint support, then $R_{a,\eps}^{t_0}$ satisfies \eqref{H1Rt}.

Let now $X$ solve between $t_0$ and $t_0+2\,\eps^\alpha$,
\[
\partial_t X(t,t_0,x)=R_{a,\eps}^{t_0}(t,X(t,t_0,x)),\quad X(t=t_0,t_0,x)=x.
\]
%
At time $t_2=t_0+2\,\eps^\alpha$, the position $X(t_2)$ is hence the result of a random transform in $\{\eps/2<|x-a|<\eps\}$ followed by random transforms in each $B(z_i,r_i\,\eps)$. 

Denote
\[\begin{split}
 &A_{a,\eps}=\{\eps\,(1-\eps^\beta)<|a-x|<\eps\}\cup\left(\bigcup_{i=1}^k \{r_i\,\eps\,(1-\eps^\beta)<|z_i-x|<r_i\,\eps\}\right),\\
&\tilde B_{a,\eps}=\bigcup_{i=1}^k B(z_i,r_i\,\eps)\setminus A_{a,\eps}.\\
\end{split}\]
First of all remark that $|X-a|$ is constant through the flow of $R^{ann}_{a,\eps,1}$ and that $X$ is moved by at most $\eta\,\eps$ through the flow of $R_{z_i,\eps,r_i}$. This implies point $i$ in Lemma \ref{lemmamixing}.

Point $ii$ of Lemma \ref{lemmamixing} is straightforward as $R_{a,\eps}^{t_0}$ is supported in $B(a,2\,\eps)$ and by point $iv$ of the definition of the spheres.
Point $iv$ of Lemma \ref{lemmamixing} is similarly straightforward. 

Points $iii$, $v$ and $vi$ now follow. If $|X(t,t_0,x)-X(t,t_0,y)|\leq \eps^{1+\beta}$ for $t_0\leq t\leq t_2$ then 
they are affected by the same transform, because it is not possible for them to belong to two different balls $B(z_i,r_i,\eps)$. Consequently if they do not belong to $A_{a,\eps}$ at either $t_0$ or $t_1$ then $|X(t,t_0,x)-X(t,t_0,y)|$ is constant along the evolution.

Thus if $x,\;y\in O_{a,\eps}\setminus A_{a,\eps}$, then with probability at least $1-C\,\eps^\beta$, so are $X(t_2,t_0,x)$ and $X(t_2,t_0,y)$ which first in turn implies that $|X(t_2,t_0,x)-X(t_2,t_0,y)|=|x-y|\leq \eps^{1+\beta}$, that is point $v$ of the Lemma. 

Either $X(t_2,t_0,x)\not\in \tilde B_{a,\eps}$ implying that $\eps/2\leq |X(t_2,t_0,x)-a|\leq \eps$. Or $X(t_2,t_0,x)$ is obtained from $X(t_1,t_0,x)$ through a random rotation of center $z_i$. In that case $X(t_2,t_0,y)$ is obtained through the same random rotation. Point $ii$ of the definition of the spheres and the fact that $|X(t_2,t_0,x)-X(t_2,t_0,y)|\leq \eps^{1+\beta}$ implies that with probability at least $1/2$, both $X(t_2,t_0,x)$ and $X(t_2,t_0,y)$ belong to the annulus $\{\eps/2\leq |x-a|\leq \eps$. That gives point $iii$ of the Lemma.

Finally point $iii$ of the definition of the spheres implies that if $x,\;y\in O_{a,\eps}\setminus A_{a,\eps}$, then with probability at least $1-\eta-C\,\eps^\beta$, they will belong to one of the balls $B(z_i,\eps\,r_i)\setminus A_{a,\eps}$ and hence that $X(t_2,t_0,x)-X(t_2,t_0,y)$ is obtained through a rotation of angle uniformly distributed in $[0,\ 2\,\pi]$. This is point $vi$ and concludes the proof of Lemma \ref{lemmamixing}.

\end{document}